\documentclass{sig-alternate}
 \usepackage{amsmath}
 \usepackage{amsfonts}
 \usepackage{amssymb}

\def\BR{{\hbox{\bf R}}}
\def\R{{\hbox{\bf R}}}

\def\Q{{\hbox{\bf Q}}}

\def\P{{\hbox{\bf P}}}
\def\E{{\hbox{\bf E}}}

\font \roman = cmr10 at 10 true pt

\def\be#1{ \begin{equation}\label{#1} }
\def\bas{\begin{align*}}
\def\eas{\end{align*}}
\def\bi{\begin{itemize}}
\def\ei{\end{itemize}}

\def\Z{{\hbox{\bf Z}}}

\def\eps{\varepsilon}

\def\emph#1{{\it #1}}
\def\textbf#1{{\bf #1}}

\def\Bv{{\mathbf v}}

\def\BI{{\mathbf I}}

\def\BZ{{\mathbf Z}}

\def\tt{{\tilde t}}

\def\ep{{\epsilon}}

\def\dd{\partial}

\def \small {\scriptsize}

  \newtheorem{theorem}[subsection]{Theorem}

  \newtheorem{lemma}[subsection]{Lemma}
  \newtheorem{corollary}[subsection]{Corollary}
  \newtheorem{remark}[subsection]{Remark}

  \newtheorem{definition}[subsection]{Definition}
\include{psfig}

\begin{document}
\conferenceinfo{STOC'07,} {June 11--13, 2007, San Diego, California, USA.}
\CopyrightYear{2007}
\crdata{978-1-59593-631-8/07/0006}

\title{The condition number of a randomly perturbed matrix }

\numberofauthors{2}
\author{
\alignauthor Terence Tao\titlenote{T. Tao is supported by a grant from the Macarthur Foundation.}\\
\affaddr{Department of Mathematics, UCLA}\\
\affaddr{Los Angeles CA 90095-1555, USA.}\\
\email{tao@math.ucla.edu}
\alignauthor Van Vu\titlenote{V. Vu is an A. Sloan  Fellow and is supported by an NSF Career Grant.}\\
\affaddr{Department of Mathematics, Rutgers}\\
\affaddr{Piscataway, NJ 08854, USA.}\\
\email{vanvu@math.rutgers.edu}
}
\maketitle

\begin{abstract} Let $M$ be an arbitrary $n$ by $n$ matrix. We study the condition
number  a random perturbation $M+N_n$ of $M$, where $N_n$ is a
random matrix. It is shown that, under very general conditions on
$M$ and $M_n$, the condition number  of $M+N_n$ is polynomial in $n$
with very high probability. The main novelty here is that we allow
$N_n$ to have discrete distribution.
\end{abstract}

\section{Introduction}
\subsection {The condition number}
Let  $M$ be an $n \times n$ matrix,
$$\sigma_1 (M):= \sup_{x \in \R^n, \|x\|=1} \|Mx\|$$

\noindent is the largest singular value of $M$ (this parameter is
also often called the operator norm of $M$).

If $M$ is invertible, the \emph{condition number} $\kappa(M)$ is
defined as
$$\kappa(M):= \sigma_1(M)  \sigma_1 (M^{-1}). $$

The condition number  plays a
 crucial role in numerical linear algebra. The
 accuracy and stability of most algorithms used to solve the
 equation $Mx=b$ depend on $\kappa(M)$. The
 exact solution $x= M^{-1} b$, in theory, can be computed quickly (by
 Gaussian elimination, say). However, in practice computers can only present a
 finite subset of real numbers and this leads to two
 difficulties. The represented numbers cannot be arbitrary large
 of small, and there are gaps between them. A quantity which is frequently used in numerical
 analysis is  $\ep_{\hbox{machine}}$ which is half of the distance
 from $1$ to the nearest represented number. A fundamental result
 in numerical analysis \cite{BT}
 asserts that if one denotes by $\tilde x$ the result computed by
 computers, then the relative error $\frac{ \| \tilde x - x \|
 }{\|x\|}$ satisfies

 $$ \frac{ \| \tilde x - x \| }{\|x\|} = O\big( \ep_{\hbox{machine}}
 \kappa(M) \big) $$

We call $M$ {\it well conditioned} if $\kappa (M)$ is small. For
quantitative purposes, we say that an $n$ by $n$ matrix $M$ is
{\it well conditioned } if its condition number is polynomially
bounded in $n$ ($\kappa(M) \le n^C$ for some constant $C$
independent of $n$).

 In the whole paper, we think that $n$ is large and the asymptotic
 notation is used under the assumption that $n \rightarrow
 \infty$.
\subsection {Effect of noise}
 An important issue in the theory of computing is noise,
as almost all  computational processes are
 effected by it.  By the word noise, we would like to  represent all kinds of
 errors occurring in a process, due to both humans and machines, including
 errors  in measuring,  errors caused by truncations,
errors committed in transmitting  and inputting the data, etc.

It happens frequently  that while we are interested in a solving a
certain
 equation, because of the noise
 the computer actually ends
 up with solving a slightly perturbed version of it.
Our work is motivated by the following phenomenon, proposed by
Spielman and Teng \cite{ST} \vskip2mm {\bf P1:} {\it For every
input instance it is unlikely that a slight random perturbation of
that instance has large condition number. } \vskip2mm

If the input is a matrix, we can reformulate this in a more
quantitative way as follows

\vskip2mm {\bf P2:} {\it Let $M$ be an arbitrary $n$ by $n$ matrix
and $N_n$ a random
 $n$ by $n$ matrix.
 Then with high probability $M+N_n$ is well conditioned.}

\vskip2mm The crucial point here is that $M$ itself may have large
condition number. The above phenomenon gives an explanation to the
fact (which has been observed numerically for some time--see
\cite{SST})   that one rarely encounters ill-conditioned matrices
in practice. This is also the core of  Spielman-Teng smooth
analysis which we will discuss in more details in   Section
\ref{smooth}.

The goal of this paper is  to show that under very general
assumptions on $M$ and $N_n$, $M+N_n$ indeed has small condition
number with overwhelming probability. The main novelty here is
that we allow the random matrix $N_n$ to have {\it discrete }
distribution. This  is  a natural assumption for random variables
involved in digital processes. On the other hand, very little has
been known, prior to this paper, about this case. Random discrete
matrices are indeed much more difficult to analyze than their
continuous counterparts  and our analysis  is significantly
different from those used earlier for the continuous models. In
particular, it relies heavily on a new development in additive
combinatorics, the so-called Inverse Littlewood-Offord theory (see
Section \ref{tools}).

\subsection {A necessary assumption}

Suppose that we would like to show that $M+ N_n$ is well
conditioned. This requires to bound both $\|M + N_n \|$ and $\|
(M+  N_n)^{-1} \|$ by a polynomial in $n$. Let us look at the
first norm. By the triangle inequality

$$\| M\| - \|N_n \| \le \|M + N_n \| \le \|M \| + \|N_n \| .$$

In most models for random matrices,   $\| N_n \| $ is $O(\sqrt n)$
with very high probability. Thus $\|M + N_n \|$ is often dominated
by $\| M \|$. So in order to make $\kappa (M + N_n) = n^{O(1)}$,
it is natural to assume that $\| M \| = n^{O(1)}$. In fact, as
$$ \| M \|^2 = \sigma_1^2 \le \sum_{ij} m_{ij}^2 = \sum_{i=1}^n
\sigma_i ^2 \le n \sigma_1^2 = n \|M \|^2, $$ \noindent  where
$m_{ij}$ are the entries of $M$, this assumption is equivalent to
saying that all entries of $M$ are polynomially bounded. We will
make this assumption about $M$ in the rest of the paper. The main
task now is to bound the second norm, $\| (M+N_n)^{-1} \|$, from
above.

\section {The results}
\subsection{Continuous  noise}
The case when entries of  $N_n$ are i.i.d Gaussian random
variables (with mean zero and variance one) has been studied by
various authors \cite{Ede, SST}. In particular, Sankar,
Spielman and Teng \cite{SST} proved
\begin{theorem} \label{theo:conditionST} Let $M$ be an arbitrary $n$ by $n$ matrix.
 Then for any $ x >0$,
$$\P( \|(M+ N_n)^{-1} \|  \ge x) = O(\frac{\sqrt n}{x} ) . $$
\end{theorem}

It is well known that there are positive constants $c_1$ and $c_2$
such that $\P( \|N_n \| \ge c_1 \sqrt n) \le \exp (-c_2 n)$.
\begin{corollary} \label{cor:conditionST} Let $B> C +3/2$ be positive constants.
 Let $M$ be an arbitrary $n$ by $n$
matrix whose entries have absolute value at most $n^C$. Then

$$\P( \kappa (M+ N_n)  \ge n^{B}) = O(n^{-B+C+3/2} ) . $$
\end{corollary}
\begin{proof} By the assumption on $M$ and the fact about $\|N_n \|$,
$\| M+ N_n \| =O(n^{C+1})$ with probability $1- \exp(-\Omega
(n))$. By Theorem \ref{theo:conditionST}, $\| (M+N_n)^{-1} \| \le
n^{B-C-1}$ with probability $O(n^{-B + C + 3/2})$. Thus the claim
follows by the union bound. \end{proof}
\subsection {Discrete noise: Bernoulli case}
Let us now consider random variables with discrete supports. By
rescaling, we can assume that their supports lie on $\BZ$ (or
$\BZ^d$ for some $d$).
 The most basic model among random discrete matrices is
 the Bernoulli matrix, whose  entries
  are i.i.d Bernoulli random
 variables (taking values $-1$ and $1$ with probability $1/2$).

Bounding the norm of the inverse of a random discrete matrix is a
difficult task, and the techniques used for the continuous case
are no longer applicable. In fact, it is already not trivial to
prove that a random Bernoulli matrix is almost surely invertible.
Efficient bounds on the norm  of the inverse of a Bernoulli random
matrix were obtained  only very recently \cite{Rud, TVsing}.

Our first result here is the discrete  analogue of Theorem
\ref{theo:conditionST}, where the Gaussian noise is replaced by
the Bernoulli noise.

\begin{theorem} \label{theorem:conditionTV} For any constants $A$
and $C$ there is a constant $B$ such that the following holds. Let
$M$ be an integer $n$ by $n$ matrix whose entries (in absolute
values) are bounded from above by $n^C$ and $N_n$ be the $n$ by
$n$ random Bernoulli matrix. Then
$$\P( \|(M+ N_n)^{-1} \|  \ge n^B) \le n^{-A} . $$
\end{theorem}

\begin{corollary} \label{cor:conditionTV} For any constants $A$
and $C$ there is a constant $B$ such that the following holds. Let
$M$ be an arbitrary $n$ by $n$ matrix whose entries (in absolute
values) are bounded from above by $n^C$ and $N_n$ be the $n$ by
$n$ random Bernoulli matrix. Then
$$\P( \kappa (M+ N_n) \|  \ge n^B) \le n^{-A} . $$
\end{corollary}

\begin{remark}
It is useful to have the right hand side be $n^{-A}$ rather than
just $o(1)$. The reason is that in certain applications (see for
instance Section \ref{smooth}), we need to show that polynomially
many matrices have, simultaneously, small condition numbers. The
bound $n^{-A}$  guarantees that  we can achieve this by a
straightforward union-bound argument.
\end{remark}

Theorem \ref{theorem:conditionTV} is a special case of a general
theorem, which, among others, asserts that the same conclusion
still holds when we replace the Bernoulli random variable by
arbitrary symmetric random discrete variables. We present this
theorem in the next subsection.

\subsection{Arbitrary discrete noise}

 {\it Notation.} For a real number $x$,
we use $e(x)$ to denote $$\exp(2\pi i x) = \cos 2\pi x + i \sin
2\pi x. $$
\begin{definition} Let $\mu \le 1/2$ and $D$ be positive constants. A
 random variable $\xi$ is
$(\mu,D)$-bounded if there is an integer $1 \le k \le D$ such that
for any $t$
$$| \E (e( \xi t)) |\le (1-\mu) + \mu \cos 2\pi k t. $$
A random vector (matrix) is $(\mu,D)$-bounded  if its coordinates
(entries)  are independent $(\mu,D)$-bounded random variables.
\end{definition}
\begin{remark} We need to assume $\mu \le 1/2$ to guarantee that
$(1-\mu) + \mu \cos 2\pi t$ is non-negative for all $t$.
\end{remark}
\begin{theorem} \label{theorem:maincor} For  any positive constants
$\mu \le 1/2$, $ A, C$ and $D$ there is a constant $B$ such that
the following holds.  Let $M$ be a fixed integer $n$ by $n$ matrix
whose entries have absolute values most $n^C$. Let $N_n$ be an $n$
by $n$ $(\mu,D)$ -bounded random matrix whose entries have
absolute values at most $n^C$ (with probability one). Then
$$\P (\sigma_n (M+N_n) \le n^{-B} ) \le n^{-A}. $$
\end{theorem}
\begin{remark} It is useful to note that the entries of $N_n$ are not
required to have the same distribution.  This allows the
possibility that the noise at a certain location has a correlation
with the corresponding entry of the original matrix $M$. For
instance, it might be natural to expect that the noise occurring
to a larger entry of $M$ have larger variance.
\end{remark}
\vskip2mm The following lemma provides a sufficient condition for
$(\mu,D)$-boundedness.
\begin{lemma} \label{lemma:epbounded} Let $\xi$ be a symmetric discrete random variable
and assume that there is a positive  integer $s$ such that
$\P(\xi=s) \ge \ep$. Then $\xi$ is $(\ep/2, 2s)$-bounded.

\end{lemma}

\begin{proof} (Proof of Lemma \ref{lemma:epbounded}) By the symmetry of $\xi$ and the
triangle inequality
$$ |\E (e(\xi t)) | = |\sum_{m =-\infty}^{\infty} \P (\xi=m) \cos
2\pi m t | \le (1-2\ep) + |2\ep \cos 2\pi st |. $$ Using the
elementary inequality $|\cos x| \le \frac{3}{4} + \frac{1}{4} \cos
2x$ with $x= 2 \pi st$, we have
$$ (1-2\ep) + |2\ep \cos 2\pi st | \le (1-\frac{\ep}{2})
+\frac{\ep}{2} \cos 4 \pi st, $$ \noindent concluding the proof.
\end{proof}
With this lemma, one can easily check that most basic variables
are $(\mu,D)$-bounded for some constants $\mu$ and $D$. Let us
list a few examples:
\begin{itemize}
\item (Bernoulli) $\xi$ is 1 or $-1$ with probability $1/2$. We
can take $\ep=1/2$ and $s=1$. \item (Lazy coin flip) $\xi=0$ with
probability $1-\alpha$ and $1$ or $-1$ with probability
$\alpha/2$. We can take $\ep= \alpha/2$ and $s=1$. \item
(Discretized Gaussian) Define $\xi$ as follows: $\P(\xi=m) = \P
(m-1/2 \le \Xi \le m +1/2)$, where $\Xi$ is  standard Gaussian. We
can take $\ep= \P (1/2 \le \Xi \le 3/2)$ and $s=1$. \item As a
generalization of the previous example, one can consider the
discretization of any symmetric random variable.

\end{itemize}

\subsection {The general result}
Now we are going to present an even more general result, which
implies Theorem \ref{theorem:maincor}. In this result, we do not
require that the entries of the random matrix be independent.
\begin{definition} \label{definition:type}
Let $ \mu \le 1/2$ and $C,K$ be positive constants. A random
vector $X$ of length $n$ is said to be of type $(\mu, C, K)$ if
\begin{itemize}
\item (boundedness) With probability one, all coordinates of $X$
are integer with  absolute value at most $n^C$. \item
(non-degeneracy) For any unit vector $y$, $\P (|X \cdot y| \le
n^{-2}) \le 1 -\mu/2$. (This means that $X$ is not concentrated
near a hyperplane.)

\item (concentration) There are positive integers $a_1, \dots,
a_n$ with $lcm (a_1, \dots, a_m) \le n^K$ such that for any vector
$v \in \BZ^n$,

\begin{equation} \label{concentration-critical} \sup _{a \in Z} \P (X \cdot v= a) \le \int_{0}^1
\prod_{i=1}^n \big((1-\mu) + \mu \cos 2\pi a_i v_i t \big) \dd
t,\end{equation}

where $lcm (a_1 \dots, a_m)$ (least common multiple) is the
smallest positive integer divisible by all $a_i$.

\end{itemize}
\end{definition}
\begin{remark} Here and later, one should not take the absolute constants such as
$-2$ and $2$ too seriously.  We make no attempt to optimize these
constants.
 The first two conditions in the definition are quite intuitive. The
third and critical condition comes from Fourier analysis and the
reader will have a better understanding of it after reading the
next section.
\end {remark}

\begin{definition} \label{definition:ind}  A collection of $n$  random vectors $Y_1,
\dots, Y_n$ in $\R^n$ is strongly linearly independent if for any
non-zero vector $y \in \R^n$ and any $ 1\le i \le n$,
$$\P (Y_1, \dots, Y_n \,\, \hbox{independent}| Y_i=y) \le  \exp(-\Omega
(n)). $$
\end{definition}

 \begin{theorem} \label{theorem:main} (Main Theorem) For every
positive constants $\mu \le 1/2,   A, C, K$ there is a positive
constant $B$ such that the following holds.  Let $M_n$ be a random
matrix with the following two properties
\begin{itemize}
\item  The   row vectors of $M_n$ are independent random vectors
of type $(\mu,  C,K)$. \item  The column vectors of $M_n$ are
strongly linearly independent.
\end{itemize}
Then
$$\P( \sigma_n (M_n) \le n^{-B}) \le n^{-A}. $$
\end{theorem}
\begin{remark} \label{remark:frozen}
Actually in the concentration property, one can omit a few
coordinates in the product. To be more precise, we can make the
following weaker assumption:

\begin{itemize}

\item There is a subset $E$ of $\{1, \dots, n \}$ of at most
$n^{.99}$ elements and  positive integers $a_i$, $i \in \{1\dots n
\} \backslash E$ with $lcm$ at most $n^K$ such that for any vector
$v \in \BZ^n$,

\begin{equation} \label{concentration-critical1} \sup _{a \in Z} \P (X \cdot v= a) \le \int_{0}^1
\prod_{i \in \{1\dots n \} \backslash E}  \big((1-\mu) + \mu \cos
2\pi a_i v_i t \big) \dd t,\end{equation}

\end{itemize}

Remark that we do not require any control on the coordinates in
$E$. This allows us  to handle, for instance,   the case when
there are frozen entries which are not effected by noise. (In this
case we simply put these coordinates in $E$.)  This situation does
occur in practice. In particular, a zero entry is often
noise-free.
\end{remark}
\section{Proof of Theorem 2.11}
In order to derive Theorem \ref{theorem:maincor} from Theorem
\ref{theorem:main}, we first need to verify that the matrix in
Theorem \ref{theorem:maincor} is of type $(\mu, C,K)$ for some
constants $\mu, C$ and $K$. This will be done in the first two
subsections. Next, we need to verify the strong linear
independence. This will be done in the last subsection.
\subsection{Checking the concentration property}
In this subsection, we verify the concentration property in the
definition of $(\mu, C,K)$-type. This is based on the following
lemma.
\begin{lemma} \label{lemma:noncon}  Let $Z$ be an arbitrary
integer vector and  $X$ be a random $(\mu,D)$-bounded vector, both
 of length $n$.
 Then there exist positive integers $a_1, \dots, a_n$ at most
$D$ such that   for any vector $v \in \BZ^n$

$$\sup _{a \in \BZ} \P ((Z+X) \cdot v = a) \le \int_{0}^1 \prod_{i=1}^n
\big((1-\mu) + \mu \cos 2\pi a_i  v_i t \big) \dd t. $$
  \end{lemma}
\begin{proof} As $a$ can take any value, it suffices to prove the
statement for $Z=0$.  For an integer $x$, the indicator
$\BI_{x=0}$ of the event $x=0$ can be expressed, using Fourier
analysis, as
$$\BI_{x=0} = \int_0^1 e ( xt ) \dd t. $$
Let $\xi_i$, $1\le i \le n$ be the coordinates of $X$. The event
$X \cdot v=a$ can be rewritten as $\sum_{i=1}^n \xi_i v_i -a =0$.
Thus
$$\P (X \cdot v=a) = \E (\BI_{\sum_{i=1}^n \xi_i v_i -a=0}) = \E \Big(\int_0^1
e(\sum_{i=1}^n \xi_i v_i -a) t ) \dd t \Big). $$

As the $\xi_i$ are independent, the last expectation is equal to
$$\int_0^1 \exp(-2\pi a t) \prod_{i=1}^n \E e( \xi_i v_i t)
\dd t \le \int_0^1  \prod_{i=1}^n |\E (e( \xi_i v_i t))| \dd t.
$$
As $\xi_i$ is $(\mu,D)$-bounded, there is a positive integer $a_i
\le D$ such that
$$ |\E (\exp (2\pi i \xi_i v_i t)| \le (1-\mu) + \mu \cos 2\pi a_i
v_i t, $$ \noindent completing the proof. \end{proof}

\subsection{Checking the non-degeneracy property}
Let $y$ be a unit vector in $\BR^n$ and $X$ be a random
$(\mu,D)$-bounded vector of length $n$ and $Z$ be an arbitrary
integer vector of length $n$. We want to show that
$$\P (| (Z+X) \cdot y | \le n^{-2} ) \le 1 -\mu/2. $$


 If $(Z+X) \cdot y$
has absolute value at most $n^{-2}$, then $X \cdot  ny$ has
absolute value at most $n^{-1}$. As $y$ is an unit vector, one of
the coordinate of $ny$ has absolute value larger than $1$. Assume,
without loss of generality, that the first coordinate $y_1$ of
$ny$ is such large. Recall that $X= (\xi_1, \dots, \xi_n)$ where
the $\xi_i$ are independent $(\mu,D)$-bounded random variables.
Condition on $\xi_2, \dots, \xi_n$, it suffices to show that for
any interval $I$ of length $2n^{-1}$
$$\P (\xi_1 v_1 \in I ) \le 1- \mu/2. $$
\noindent But since $\xi$ take only integer values and $|y_1| \ge
1$,  the values of $\xi_1 y_1$ would be at least one apart.
Assume, for a contradiction, that  $\P (\xi_1 v_1 \in I ) >
1-\mu/2$. This would imply that there is a number $s$ such that
$\P(\xi_1=s) > 1-\mu/2$. Then by the triangle inequality
$$|\E (e(\xi_1 t)) | \ge |e(st)| (1-\mu/2) - \mu/2 \ge 1-
\mu, $$ for any $t$. On the other hand, as $\xi_1$ is
$(\mu,D)$-bounded
$$ |\E (e(\xi_1 t)) | \le (1-\mu) +\mu \cos 2\pi a_1 t $$ for some
$a_1 \le D$. Taking $t$ such that $ \cos 2\pi a_1 t = -1$, we
obtain  a contradiction and conclude the proof.

\subsection{Checking the strong linear independence}
The strong linear independence of the column vectors of a  random
$(\mu,D)$-bounded matrix is a consequence of the following
theorem, which can be proved by refining the proof of
\cite[Theorem 1.6]{TV1}.

\begin{theorem} \label{theo:KKS} Let $\mu \leq 1/2$ and $D,l$ be
positive constants.
 Then there is a positive
constant $\eps=\eps(\mu, D, l)$ such that the following holds. For
any set $Y$ of  $l$ independent vectors from $\R^n$ and $n-l$
independent random $(\mu,D)$-bounded vectors of length $n$, the
probability that they  are linearly dependent is at most $
(1-\eps)^n$.
\end{theorem}
\begin{remark} This theorem is a generalization  of
a well known theorem of Kahn, Koml\'os and Szemer\'edi \cite{KKS}
which asserts that the probability that  a random Bernoulli matrix
is singular is exponentially small. To see this, recall that a
random Bernoulli vector is $(1/4,2)$-bounded and in Theorem
\ref{theo:KKS} take $l=1$ and fix $y$ be the all one vector.
\end{remark}

\section{  Smooth complexity with discrete noise} \label{smooth}
Running times of algorithms are frequently estimated by worst-case
analysis. But in practice, it has been observed that many
algorithms perform significantly better than the estimates
obtained from the worst-case analysis. Few years ago, Spielman and
Teng \cite{ST, ST1} came up with an ingenuous explanation for this
fact. The rough idea behind their argument is as follows. Even if
the input $I$ is the worst-case one (which, in theory, would
require a long running time), because of the noise, the computer
actually works on some slightly randomly perturbed version of $I$.
Next, one would  show that the running time on a slightly randomly
perturbed input, with high probability, is much smaller than the
worst-case one.  The smooth complexity of an algorithm is the
maximum over its input of the expected running time of the
algorithm under slight perturbations of that input. The puzzling
question here is, of course: why the perturbed input is typically
better than the original (worst-case) one ? In some sense, the
"magic" lies in the Phenomenon {\bf P1}. The random noise
guarantees that the condition number of the perturbed input is
small (so the perturbed input is likely to be well conditioned),
no matter how ill conditioned the original input may be. The bound
on the condition number then can be used to derive a bound on the
running time of the algorithm.

In their works \cite{ST, ST1, SST}, Spielman and Teng (and
coauthors) assumed Gaussian noise (or more generally continuous
noise). Theorem \ref{theo:conditionST} played a significant role
in their proofs.

  An important (and largely open) problem is to obtain smooth complexity bounds
   when the noise is discrete. (We would like to thank Spielman for
communicating this problem.) In fact, it is not clear how
computers would compute with  Gaussian (and other continuous)
distributions without discretizing them. This problem seems to
pose a considerable mathematical challenge. Naturally, the first
step would be to obtain estimates for the condition number with
discrete noise. This step has now been accomplished in this paper.
However, these estimates themselves are not always  sufficient. To
be more specific, the situation looks as follows:

\begin{itemize}
\item There are problems where an efficient bound on the
condition number leads directly to an efficient complexity bound.
In such a  situation, we obtain a smooth complexity bound with
discrete noise in the obvious manner. This seems to be the case,
e.g., with the problems involving the Gaussian Elimination in
\cite{SST}. In the proofs in \cite{SST}, the critical fact was
that all $n-1$ minors of a random perturbed matrix are all well
conditioned, with high probability. This can be obtained using our
results combined with the union bound (see the remark after
Theorem \ref{theorem:conditionTV}).

\item There are situations where beside the estimate on the
condition number, further properties of the noise is used. An
important example is the simplex method in linear programming. In
the smooth analysis of this algorithm with Gaussian noise
\cite{ST1}, the fact that  the distribution is continuous was
exploited  at several places. Thus, even with the discrete version
of the condition number estimates in hand, it is still not clear
to us  how to obtain a smooth complexity bound with discrete noise
in this problem.
\end{itemize}

\section  { Key ingredients} \label{tools}

In this section, we present our key ingredients in the proof of
Theorem \ref{theorem:main}.

\subsection{Generalized arithmetic progressions and their discretization}
\label{discrete-subsec}

One should take care to distinguish the sumset $kA$ from the dilate
$k \cdot A$, defined for any real $k$ as
$$k \cdot A := \{ ka| a\in A \}. $$

Let $P$ be a GAP of integers of  rank $d$ and volume $V$. Our first
key ingredient is a theorem that shows that given any specified
scale parameter $R_0$, one can ``discretize'' $P$ near the scale
$R_0$. More precisely, one can cover $P$ by the sum of a coarse
progression and a small progression, where the diameter of the small
progression is much smaller (by an arbitrarily specified factor of
$S$) than the spacing of the coarse progression, and that both of
these quantities are close to $R_0$ (up to a bounded power of $SV$).
\begin{theorem}[Discretization]\label{discrete-thm}\cite{TVsing}  For every
constant $d$ there is a constant $d'$ such that the following hold.
 Let $P \subset \Z$ be a symmetric generalized arithmetic
 progression of rank $d$ and volume $V$.  Let $R_0, S$ be positive integers.
Then there exists a number $R \geq 1$ and two generalized
progressions $P_{\operatorname{small}}$, $P_{\operatorname{sparse}}$
of \emph{rational} numbers with the following properties.
\begin{itemize}
\item(Scale) We have $R \le (SV)^{d'} R_0$. \item(Smallness)
$P_{\operatorname{small}}$ has rank at most $d$, volume at most $V$,
and takes values in $[-R/S, R/S]$. \item(Sparseness)
$P_{\operatorname{sparse}}$ has rank at most $d$, volume at most
$V$, and any two distinct elements of $SP_{\operatorname{sparse}}$
are separated by at least $RS$. \item(Covering) We have $P \subseteq
P_{\operatorname{small}} + P_{\operatorname{sparse}}$.
\end{itemize}
\end{theorem}

\subsection{Inverse Littlewood-Offord theorem}
Our second key ingredient is a theorem which characterizes all sets
$\Bv= \{v_1, \dots, v_n \}$ such that $\int_{0}^1 \prod_{i=1}^n
\big((1-\mu) + \mu \cos 2\pi  v_i t \big) \dd t$ is large. This
theorem is a refinement of \cite[Theorem 2.5]{TVsing} (see Remark
2.8 from this paper)  and will enable us to exploit the
non-concentration property from Definition \ref{definition:type} in
a critical way.

\begin{theorem} \label{theorem:TVsing} Let $0 < \mu \leq 1$ and $A, \alpha > 0$ be arbitrary.
Then there is a positive constant  $A'$ such that the following
holds.
 Assume that $\Bv = \{v_1, \ldots, v_n\}$ is a multiset of integers satisfying
$$\int_{0}^1 \prod_{i=1}^n \big((1-\mu)+ \mu \cos 2\pi v_i \xi
\big) \dd \xi \geq n^{-A}. $$

\noindent Then there is a GAP $Q$ of rank at most $A'$ and volume at
most $n^{A'}$ which contains all but at most $n^{\alpha}$ elements
of $\Bv$ (counting multiplicity). Furthermore, there is a integer $1
\le s \le n^{A'}$ such that $su \in \Bv$ for each generator $u$ of
$\Q$. \end{theorem}

With the two key tools in hand, we are now ready to prove  Theorem
\ref{theorem:main}.
\section{ Proof of Theorem 2.18}

Let $B
> 10$ be a large number (depending on the type of $M_n$) to be chosen later.
If  $\sigma_n M_n < n^{-B} $ then there exists a unit vector $v$
such that
$$ \| M _n v \| < n^{-B}.$$
By rounding each coordinate $v$ to the nearest multiple of
$n^{-B-2}$, we can find a vector $\tilde v \in n^{-B-2} \cdot \Z^n$
of magnitude $0.9 \le \|\tilde v \| \le 1.1$ such that
$$ \| M_n \tilde v \| \le 2n^{-B}.$$

Writing $w := n^{B+2} \tilde v$, we thus can find an integer vector
$w  \in \Z^n$ of magnitude $ .9 n^{B+2} \le\|w\| \le 1.1 n^{B+2} $
such that
$$ \| M_n w \| \le 2n^2.$$
Let $\Omega$ be the set of integer vectors $w  \in \Z^n$ of
magnitude $ .9 n^{B+2} \le\|w\| \le 1.1 n^{B+2} $. It suffices to
show the probability bound
$$ \P( \hbox{there is some } w \in \Omega \hbox{ such that } \| M_n w\| \le 2n^2 )
\le n^{-A} .$$ We now partition the elements $w = (w_1,\ldots,w_n)$
of $\Omega$ into three sets:

\begin{itemize}
\item  We say  that $w$ is \emph{rich} if
$$ \sup_{a \in \BZ, 1\le i \le n} \P(X_i \cdot w =a)  \geq n^{-A-4}, $$
\noindent where $X_i$ are the  row vectors of $M_n$.  Otherwise we
say that $w$ is \emph{poor}. Let $\Omega_1$ be the set of poor
$w$'s.

\item A rich $w$ is \emph{singular} $w$ if fewer than $n^{0.2}$ of
its coordinates have absolute value $n^{B/2}$ or greater. Let
$\Omega_2$ be the set of rich and singular $w$'s.

\item A rich $w$ is \emph{non-singular} $w$, if at least $n^{0.2}$
of its coordinates have absolute value $n^{B/2}$ or greater. Let
$\Omega_3$ be the set of rich and non-singular $w$'s.
\end{itemize}
\begin{remark} Again one should not take the absolute constants
$-4, 1/2$ and $.2$ too seriously. \end{remark} \noindent The desired
estimate follows directly  from the following lemmas and the union
bound.

\begin{lemma}[Estimate for poor $w$]\label{lemma:Omega1}
$$\P( \hbox{there is some} \,\, w \in \Omega_1 \,\,\hbox{such
that} \,\,\|M _n w\| \le 2n^2 ) = o( n^{-A} ). $$
\end{lemma}

\begin{lemma}[Estimate for rich singular $w$]\label{lemma:Omega2}
$$\P( \hbox{there is some} \,\, w \in \Omega_2 \,\,\hbox{such
that} \,\,\|M _n w\| \le 2n^2 ) =o(n^{-A} ).
$$
\end{lemma}

\begin{lemma}[Estimate for rich non-singular $w$]\label{lemma:Omega3}
$$\P( \hbox{there is some} \,\, w \in \Omega_3 \,\,\hbox{such
that} \,\,\|M^{\ast} _n w\| \le 2n^2 ) =o(n^{-A} ).
$$
\end{lemma}

The proofs of Lemmas \ref{lemma:Omega1} and \ref{lemma:Omega2} are
relatively simple and rely on well-known methods.  The proof of
Lemma \ref{lemma:Omega3}, which is essentially the heart of the
matter, is more difficult and requires the tools provided in Section
\ref{tools}.
\section {Proof of Lemma 7.2}

We use a conditioning argument, following \cite{Rud}. (An argument
of the same spirit was used by Koml\'os to prove the bound
$O(n^{-1/2})$ for the singularity problem \cite{Bol}.) Let $M$ be a
matrix such that there is $w \in \Omega_1$ satisfying $\|Mw\| \le
2n^2$. Since $M^{-1}$ and its transpose have the same spectral norm,
there is a vector $w'$ which has the same norm as $w$ such that
$\|w'M\| \le 2n^2$. Let $u=w'M$ and $X_i$ be the row vectors of $M$.
Then
$$u =\sum_{i=1}^n w_i' X_i $$
\noindent where $w_i'$ are the coordinates of $w'$. Now consider
$M=M_n$. By paying a factor of $n$ in the probability (whenever this
phrase is used, keep in mind that we will use the union bound to
conclude the proof), we can assume that $w'_n$ has the largest
absolute value among the $w_i'$. We expose the first $n-1$ rows
$X_1, \dots, X_{n-1}$ of $M_n$. If there is $w \in \Omega_1$
satisfying $\|Mw\| \le 2n^2$, then there is a vector $y \in
\Omega_1$, depending only on the first $n-1$ rows such that
$$(\sum_{i=1}^{n-1} (X_i \cdot y)^2 )^{1/2} \le 2n^2. $$
\noindent  We can write $X_n$ as

$$X_n =\frac{1}{w_n'} (u- \sum_{i=1}^{n-1} w_i' X_i ). $$

\noindent Thus,
$$|X_n \cdot y| = \frac{1}{|w_n'|} |u \cdot y  - \sum_{i=1}^{n-1} w_i' X_i \cdot y|.
$$

\noindent The right hand side, by the triangle inequality, is at
most
$$\frac{1}{|w_n'|} (|u||y| + \|w'\| (\sum_{i=1}^{n-1} (X_i
\cdot y)^2 )^{1/2}). $$

\noindent By assumption $|w_n'| \ge n^{-1/2} |w'|$. Furthermore, as
$|u| \le 2n^2$, $|u||y| \le 2n^2 |y| \le 3n^2 |w'|$ as $|w'|=|w|$
and both $y$ and $w$ belong to $\Omega_1$. (Any two vectors in
$\Omega_1$ has roughly the same length.) Finally $(\sum_{i=1}^{n-1}
(X_i \cdot y)^2 )^{1/2} \le 2n^2$. Putting all these together, we
have
$$|X_n \cdot y| \le 5n^{5/2}. $$
\noindent Recall that  both $X_n$ and $y$ are integer vectors, so
$X_n \cdot y$ is an integer. The probability that $|X_n \cdot y| \le
5n^{5/2}$ is at most
$$(10n^{5/2}+1) \sup_{a \in \BZ} \P (X_n \cdot y \in I). $$
 On the other
hand, $y$ is poor, so by definition $\sup_{a \in \BZ} \P(X_n \cdot y
=a) \le n^{-A-4}$. Thus, it follows that
\begin{align*}\P( \hbox{there is some} \,\, w \in& \Omega_1 \,\,\hbox{such
that} \,\,\|M _n w\| \le 2n^2 ) \le\\
&\le n^{-A-4} (10 n^{5/2}+1) n = o( n^{-A} ),
\end{align*}
\noindent where the extra factor $n$ comes from the assumption that
$w_n'$ has the largest absolute value. This completes the proof.

\section { Proof of Lemma 7.3}
\label{section:Omega12}
 We use an argument from \cite{Lit}.  The key point will be that
the set $\Omega_2$ of rich non-singular vectors has sufficiently low
entropy that one can proceed using the union bound. A set $N$ of
vectors on the $n$-dimensional unit sphere $S_{n-1}$ is said to be
an \emph{$\ep$-net} if for any $x \in S_{n-1}$, there is $y\in N$
such that $\|x-y\| \le \ep$.  A standard greedy argument shows
\begin{lemma} \label{lemma:net} For any $n$ and $\epsilon \le 1$, there exists an $\epsilon$-net
of cardinality at most $O(1/\eps)^n$.
\end{lemma}
We need another lemma, showing that for any unit vector $y$, very
likely  $\| M_n y \|$ is polynomially large.
\begin{lemma} \label{lemma:large}
For any unit vector $y$
$$\P (\| M_n y \| \le n^{-2} ) = \exp(-\Omega (n)). $$
\end{lemma}

\begin{proof} If $\|M_ny \| \le n^{-2}$, then
$|X_i \cdot y| \le n^{-2} $ for all index $1\le i \le n$. However,
by the assumption of the theorem, for any fixed $i$ , the
probability that $|X_i \cdot y| \le n^{-2} $ is at most $1-\mu/2$.
Thus,
$$\P (\| M_n y \| \le n^{-2} ) \le  (1- \mu/2)^{n}
=\exp(-\Omega (n)) $$ concluding the proof.
\end{proof}
For a vector $w \in \Omega_2$, let $w'$ be its normalization
$w':w/\|w\|$. Thus, $w'$ is an unit vector with at most $n^{0.2}$
coordinates with absolute values larger or equal $n^{-B/2}$. By
choosing $B \ge 2C +20$, we can assume that $w'$ belong to
$\Omega_2'$, the collection of unit vectors at most $n^{0.2}$
coordinates with absolute values larger or equal $n^{-C-10}$. If $\|
Mw\| \le 2n^2$ for some $w\in \Omega_2$, then $\|Mw'\| \le 3 n^{-B}$
, as $\|w\| \ge .9 n^{B+2}$. Thus, it suffices to give an
exponential bound on the event that there is $w'\in \Omega_2'$ such
that $\| M _n w'\| \le 3n^{-B}$. By paying a factor of $\binom{n}{
{n^{0.2}} }= \exp(o(n))$ in probability, we can assume that the
large coordinates (with absolute value at least $n^{-C-10} $) are
among the first $l := n^{0.2}$ coordinates. Consider an
$n^{-C-5}$-net $N$ in $S_{l-1}$. For each vector $y \in N$, let $y'$
be the $n$-dimensional vector obtained from $y$ by letting the last
$n-l$ coordinates be zeros, and let $N'$ be the set of all such
vectors obtained. These vectors have magnitude between $0.9$ and
$1.1$, and from Lemma \ref{lemma:net} we have $|N'| \leq
O(n^3)^{l}$. Now consider a rich singular vector $w' \in \Omega_2$
and let $w^{''}$ be the $l$-dimensional vector formed by the first
$l$ coordinates of this vector. As the remaining coordinates are
small $\|w^{''}\| =1 +O(n^{-C-9})$. There is a vector $y \in N$ such
that
$$\|y - w^{''} \| \le n^{-C-5}+ O(n^{-C-9}) . $$
\noindent It follows that there is a vector $y' \in N'$ such that
$$\|y' - w'\| \le  n^{-C-5}+ O(n^{-C-9}) \le 2n^{-C-5}. $$
\noindent If $M$ has norm at most  $n^{C+1}$, then
$$\|M w'\| \ge \|M y'\| - 2n^{-C-5} n^{C+1} = \|My'\| - 2n^{-4}. $$
\noindent It follows that if $\|M w'\| \le 3n^{-B}$ for some $B \ge
2$, then $\| My'\| \le 5n^{-4}$. Now take $M=M _n$. For each fixed
$y'$, the probability that $\| M_n y'\| \le 5n^{-4} \le n^{-2} $ is
at most $\exp(-\Omega(n))$, by Lemma \ref{lemma:large}. Furthermore,
the number of $y'$ is subexponential (at most $O(n^{C+3})^{l}
O(n)^{3 n^{.2}} = \exp(o(n))$). The claim follows by the union
bound.

\section{Proof of Lemma 7.4}
\label{section:Omega3}
 This is the most difficult part of the proof, where we will need
all  the tools provided in Section \ref{tools}. Informally, the
strategy is to use the inverse Littlewood-Offord theorem  to place
the integers $w_1,\ldots,w_n$ in a progression, which we then
discretize using Theorem \ref{discrete-thm}.  This allows us to
replace the event $\|M _n w\| \le 2n^2$ by some dependence event
involving the columns of $M_n$, whose probability is very small by
the strong linear independence assumption of the theorem.

We now turn to the details. By the inverse theorem and the
non-concentration property from Definition \ref{definition:type},
there is a constant  $A'$ such that for each $w \in \Omega_3$ there
exists a symmetric GAP $Q$ of integers of rank at most $d$ and
volume at most $n^{A'}$ and non-zero integers $a_1, \dots, a_n$ with
least common multiple at most $n^K$ such that $Q$ contains all but
$\lfloor n^{0.1}\rfloor$ of the integers $a_1w_1,\ldots, a_nw_n$.
Furthermore,  the generators of $Q$ are of the form $a_iw_i/s$ for
some  $1 \leq s \leq n^{A'}$. Notice that if $a_i w_i \in Q$ then
$w_i \in Q':= \{ x/a | x \in Q, a \in \BZ, a \neq 0, |a| \le n^{K}
\}$. Using the description of $Q$ and the fact that $w_1, \dots,
w_n$ and $a_1, \dots, a_n$ are polynomially bounded in $n$, one can
see that the total number of possible $Q$ is $n^{O(1)} =\exp(o(n))$.
Next, by paying a factor of
$$ \binom{n}{\lfloor n^{0.1}\rfloor} \leq n^{\lfloor n^{0.1}\rfloor} = \exp(o(n))$$
\noindent we may assume that it is the last $\lfloor n^{0.1}\rfloor$
integers $a_{m+1}w_{m+1}$,\\$ \ldots, a_nw_n$ which possibly lie
outside $Q$, where we set $m := n-\lfloor n^{0.1} \rfloor$. As each
of the $w_i$ has absolute value at most $1.1 n^{B+2}$, the number of
ways to  fix these exceptional elements is at most $(2.2 n^{B+2}
)^{n^{0.1}} = \exp(o(n))$. Overall, it costs a factor only
$\exp(o(n))$ (keep in mind that we intend to use the union bound) to
fix $Q$, the positions and values of the exceptional elements of
$w$.

\noindent Notice that $M_n w = w_1 Y_1 + \dots w_n Y_n$, where $Y_i$
is the $i$th column of $M_n$. Fixing $w_{m+1}, \dots, w_n$ and set
$Y:= \sum_{i=m+1}^n w_i Y_i$. This way we can rewrite $M_n w$ as

$$ M_n w = w_1 Y_1 + \ldots + w_{m} Y_{m} + Y.$$
\noindent  For any number $y$, define $F_y$ be the event that there
exists $w_1,\dots, w_m$ in the set $ Q'$, where at least one of the
$w_i$ has absolute value larger or equal $n^{B-10}$, such that
$$ |w_1 Y_1 + \ldots + w_{m} Y_{m} + y| \le 2n^2 .$$
\noindent It suffices to prove that for any $y$
$$\P( F_y) \le \exp(-\Omega( n) ). $$
We now apply Theorem \ref{discrete-thm} to the GAP $Q$ with $R_0 :=
n^{B/3}$ and $S := n^{L}$, where $L=C+K+2$ ($C$ and $K$ are the
constants in Definition \ref{definition:type}). By choosing $B$
sufficiently large, we can guarantee that $B/3$ is considerably
larger than $L$.  Recall that the volume of $Q$ is at most $n^{A'}$,
where $A'$ is a constant depending on $A$ and $\mu$. We can find a
number $R = n^{B/3 + O_{A',L}(1)}$ and symmetric GAPs
$Q_{\operatorname{sparse}}$, $Q_{\operatorname{small}}$ of rank at
most $d' =d'(d, A')$ and volume at most $n^{A'}$ such that
\begin{itemize}
\item $ Q \subseteq Q_{\operatorname{sparse}} +
Q_{\operatorname{small}}.$ \item $ Q_{\operatorname{small}}
\subseteq [-n^{-L} R, n^{-L} R]$. \item The elements of $n^{L}
Q_{\operatorname{sparse}}$ are $n^{L} R$-separated.

\end{itemize}
Since $Q$ (and hence $n^{L} Q$) contains $a_1w_1,\ldots, a_mw_m$
(for some set $\{a_1, \dots, a_m\}$)  we can therefore write
$$ w_j = a_j^{-1} (w^{\operatorname{sparse}}_j + w^{\operatorname{small}}_j)$$
for all $1 \leq j \leq m$, where $w^{\operatorname{sparse}}_j \in
Q_{\operatorname{sparse}}$ and $w^{\operatorname{small}}_j \in
Q_{\operatorname{small}}$. In fact, this decomposition is unique.
Suppose that the event $F_y$ holds. Let  $y = (y_1,\ldots,y_n)$ and
$\eta_{i,j}$ denote the entry of $M_n$ at row $i$ and column $j$. We
have
$$ w_1 \eta_{i,1} + \ldots + w_m \eta_{i,m} = y_i + O(n^2).$$
\noindent for all $1 \leq i \leq n$.  Split the $w_j$ into sparse
and small components and estimating the small components. The
contribution coming from the small components is

$$ \sum_{j=1}^m a_j^{-1} w^{\operatorname{small}}_j
\eta_{i,j} = O(n^{-L+C+1} R) $$ \noindent since  $|\eta_{i,j}|$ are
bounded from above by $n^C$, $|w^{\operatorname{small}}_j|$ is
bounded from above by $n^{-L} R$ and $a_j$ are positive integers. By
the triangle inequality, it follows that
$$ a_1^{-1}  w^{\operatorname{sparse}}_1 \eta_{i,1} + \ldots + a_m^{-1}
w^{\operatorname{sparse}}_m \eta _{i,m} = y_i + O(n^{-L+C+1} R)$$
for all $1 \leq i \leq n$.

Set $T:= lcm(a_1, \dots, a_m)$. The previous estimate implies
$$b_1 w^{\operatorname{sparse}}_1 \eta_{i,1} + \ldots + b_m
w^{\operatorname{sparse}}_m \eta _{i,m} =T y_i + O(Tn^{-L+C+1} R)$$
\noindent where $b_i = T/a_i$. Now we use the assumption that $T \le
n^K$ from Definition \ref{definition:type}.  This assumption yields
that  $b_i \le n^K$ and the
  left-hand side lies in
  $$n^{K+1}Q_{\operatorname{sparse}} \subset
n^{K+1} Q_{\operatorname{sparse}} \subset n^{L} Q, $$ which is known
to be $n^{L} R$-separated. Furthermore,
$$ O(Tn^{-L+C+1}
R) = O(n^{K-L+C+1} R) = O(n^{R-1})$$ by the definition of $L$. Thus
there is a unique value for the right-hand side, call it $y'_i$,
which depends only on $y$ and $Q$
 such that

$$ b_1 w^{\operatorname{sparse}}_1 \eta_{i,1} + \ldots + b_m
w^{\operatorname{sparse}}_m \eta_{i,m} = y'_i. $$

\noindent The point is that we have now eliminated the $O()$ errors,
and have thus essentially converted the singular value problem to a
problem about dependence. Note also that since one of the
$w_1,\ldots,w_m$ is known to have magnitude at least $n^{B/2}$
(which will be much larger than $n^{L} R = n^{L + B/3}$  given that
we  set $B > 6L= 6(C+K+2)$), we see that at least one of the
$w_1^{\operatorname{sparse}},\ldots,w_n^{\operatorname{sparse}}$ is
non-zero.

Let $y' =(y_1', \dots, y_n')$. The equation
$$ b_1 w^{\operatorname{sparse}}_1 \eta_{i,1} + \ldots + b_m
w^{\operatorname{sparse}}_m \eta_{i,m} = y'_i $$ implies that the
first $m$ columns of $M_n$ span $y'$. For any fixed non-zero $y'$,
the probability that this happens is exponentially small by the
strong linear independence assumption. This completes the proof.

\section{ Frozen entries}

We now give an explanation to Remark \ref{remark:frozen}. This
remark is based on the fact that in the previous proof  one is
allowed to have as many as $n^{1-\ep}$ coordinates outside the set
$Q'$, for any positive constant $\ep <1$. Indeed, these extra
coordinates  contribute a factor of $\binom{ n}{ {n^{1-\ep}} }$
which is $\exp(o(n))$. This factor will be swallowed by the
exponential bound we have at the end of the proof. (In the proof we,
for convenience,  set $\ep =.9$ and have $n^{.1}$ exceptional
coordinates, but the actual value of $\ep$ plays no role.) The main
point here is that we can set aside the "frozen" coordinates even
before applying the Inverse Littlewood-Offord theorem.


\begin{thebibliography}{10}
\bibitem{BT} D. Bau and L. Trefethen, Numerical linear algebra.
Society for Industrial and Applied Mathematics (SIAM),
Philadelphia, PA, 1997.

\bibitem{Bol} B. Bollob\'as, Random graphs. Second edition,
Cambridge Studies in Advanced Mathematics, 73. Cambridge
University Press, Cambridge, 2001.


\bibitem{Ede}  A. Edelman,  Eigenvalues and condition numbers of random matrices.
{\it SIAM J. Matrix Anal. Appl.}  9 (1988), no. 4, 543--560.

\bibitem{ET} A. Edelman and B. Sutton, Tails of condition number distributions,
{\it SIAM J. Matrix Anal. Appl.}  27 (2005), no. 2, 547--560.


\bibitem{KKS}
J. Kahn, J. Koml\'os, E. Szemer\'edi, On the probability that a
random $\pm 1$ matrix is singular, \emph{J. Amer. Math. Soc.}
\textbf{8} (1995), 223--240.


\bibitem{Lit} A. Litvak, A. Pajor, M. Rudelson and N. Tomczak-Jaegermann,
 Smallest singular value of random matrices and geometry of random polytopes,
 {\it  Adv. Math. } 195 (2005), no. 2, 491--523.

\bibitem{Rud} M. Rudelson,  Invertibility of random matrices: Norm of the inverse.
     {\it submitted}.
\bibitem{SST}  A. Sankar, S. H. Teng, and D. A. Spielman,
Smoothed Analysis of the Condition Numbers and Growth Factors of
Matrices, {\it preprint}.

\bibitem{ST}  D. A. Spielman and S. H. Teng,
 Smoothed analysis of algorithms, {\it  Proceedings of the International Congress of
Mathematicians}, Vol. I (Beijing, 2002), 597--606, Higher Ed.
Press, Beijing, 2002.
\bibitem{ST1}  D. A. Spielman and S. H. Teng, Smoothed analysis of algorithms: why the simplex algorithm
usually takes polynomial time, {\it J. ACM } 51 (2004), no. 3,
385--463.
\bibitem{TV1} T. Tao and V. Vu, On random $\pm 1$ matrices: Singularity
and  Determinant,  {\it Random Structures Algorithms } 28 (2006),
no. 1,
 1--23.

\bibitem{TVsing} T. Tao and V. Vu, Inverse Littlewood-Offord theorems and the condition number of
 random discrete matrices, {\it Annals of Mathematics, to appear}.

\end{thebibliography}
\end{document}